\documentclass{amsart}
\usepackage{math,epsfig,psfrag}
\usepackage[all]{xy}
\renewcommand\S{{\widehat\C}}
\newcommand\M{{\mathbb M}}
\newcommand\J{{\mathfrak J}}
\def\pair|#1|{{\ll}#1{\gg}}

\begin{document}
\author{Laurent Bartholdi\and Rostislav I. Grigorchuk}
\title{On a group associated to $z^2-1$}
\date\today
\begin{abstract}
  We construct a group $G$ acting on a binary rooted tree; this
  discrete group mimics the monodromy action of iterates of
  $f(z)=z^2-1$ on associated coverings of the Riemann sphere.

  We then derive some algebraic properties of $G$, and describe for
  that specific example the connection between group theory, geometry
  and dynamics.

  The most striking is probably that the quotient Cayley graphs of $G$
  (aka ``Schreier graphs'') converge to the Julia set of $f$.
\end{abstract}
\maketitle

\section{Introduction}
The purpose of this paper is to introduce, in a very concrete context,
the construction of the \emph{iterated monodromy group} of a branched
covering of a Riemann surface. Our main focus will be on the covering
of the Riemann sphere given by the map $z\mapsto z^2-1$, whose group $G$
possesses many interesting algebraic and analytic properties. Our main
guideline is that these properties derive from dynamic properties of
the covering.

We will show:
\begin{thm}
  Let $G$ be the iterated monodromy group (see
  Definition~\ref{defn:img}) of $z^2-1$. Then $G$ acts on the binary
  rooted tree, and is generated by
  \[a=\pair|b,1|(1,2),\quad b=\pair|a,1|.\]
  where $\sigma$ permutes the top two branches of $G$, and
  $\pair|x,y|$ describes the automorphism acting respectively as $x$
  and $y$ on the left and right subtrees.

  \begin{enumerate}
  \item $G$ is weakly branch (Proposition~\ref{prop:wb});
  \item $G$ is a torsion-free group (Proposition~\ref{prop:tf});
  \item $G$ contains a free monoid of rank $2$, and hence has
    exponential growth (Proposition~\ref{prop:fm});
  \item $G$ does not contain any non-abelian free group
    (Proposition~\ref{prop:nofg});
  \item $G$ is not solvable, but every proper quotient of $G$ is
    (free abelian)-by-(finite $2$-group) (Proposition~\ref{prop:jns});
  \item $G$ has a presentation (Proposition~\ref{prop:lpres})
    \[G = \langle a,b|\,[[a^p,b^p],b^p],[[b^p,a^{2p}],a^{2p}]\text{
      for all $p$ a power of $2$}\rangle.\]
    The Schur multiplier $H_2(G,\Z)$ is free abelian of infinite rank.
  \item $G$ has solvable word problem (Proposition~\ref{prop:wcp});
  \item $G$'s spectrum is a Cantor set, and is the intersection of a
    line and the Julia set of a degree-$2$ polynomial mapping in
    $\R^3$ (Proposition~\ref{prop:spectrum})
    \[F(\lambda,\mu,\nu)=(\lambda^2+2\lambda\nu-2\mu^2,\lambda\nu+2\nu^2,-\mu^2).\]
  \item $G$'s limit space is the Julia set $\J$ of $z^2-1$ , and the
    Schreier graphs $G/P_n$ converge to $\J$
    (Proposition~\ref{prop:limit}).
  \end{enumerate}
\end{thm}

The group $G$ has many properties in common with the group $H$
introduced by Brunner, Sidki and
Vieira~\cite{brunner-s-v:nonsolvable}; the main difference is the
existence of free monoids in $G$, which do not seem to exist in $H$.
We also believe that most proofs for $G$ are simpler than the
analogous statements for $H$; in particular, the spectrum on $H$ has
not been computed.

\section{Covering Maps and Groups}
This section describes the general construction of a group from a
rational map of the Riemann sphere $\S$.

Let $f\in\C(z)$ be a branched self-covering of $\S$. A point $z\in\S$
is \emph{critical} if $f'(z)=0$, and is a \emph{ramification point} if
it is the $f$-image of a critical point. The \emph{postcritical set}
of $f$ is $\{f^n(z):\,n\ge1,\,z\text{ critical}\}$.
  
Assume $P$ is finite, and write $\M=\S\setminus P$. Then $f$ induces
by restriction a self-covering of $\M$.

Let $*$ be a generic point in $\M$, i.e.\ be such that the iterated
inverse $f$-images of $*$ are all distinct. If $f$ has degree $d$,
then for all $n\in\N$ there are $d^n$ points in $f^{-n}(*)$, and these
form naturally the $n$th \emph{layer} $T_n$ of a $d$-regular rooted
tree $T=\bigcup_{n\ge0}T_n$. The root of $T$ is $*$, and there is an
edge between $z$ and $f(z)$ for all $z\in T\setminus\{*\}$. In
particular, write $T_1=f^{-1}(*)$ the first layer of $T$.

Let $\gamma$ be a loop at $*$ in $\M$. Then for all $v\in T_n$ there
is a unique lift $\gamma_v$ of $\gamma$ starting at $v$ such that
$f^n(\gamma_v)=\gamma$; and furthermore the endpoint $v^\gamma$ of
$\gamma_v$ also belongs to $T_n$.
\begin{prop}
  For any such $\gamma$ the map $v\mapsto v^\gamma$ is tree
  automorphism of $T$, and depends only on the homotopy class of
  $\gamma$ in $\pi_1(\M,*)$.
\end{prop}

\begin{defn}\label{defn:img}
  The \emph{iterated monodromy group} $G_T(f)$ of $f$ is the subgroup
  of $\aut(T)$ generated by all maps $v\to v^\gamma$, as $\gamma$
  ranges over $\pi_1(\M,*)$.
\end{defn}

This definition is actually independent of the choice of $*$:
\begin{prop}
  Let $*'$ be another generic basepoint, with tree $T'$, and choose a
  path $p$ from $*$ to $*'$. Then there is an isomorphism $\phi:T\to
  T'$ such that
  \[\xymatrix{{\pi_1(\M,*)}\ar[r]^{p_*}\ar[d]_{\text{act}} &
    {\pi_1(\M,*')}\ar[d]^{\text{act}}\\
    {G_T(f)}\ar[r]_{\phi_*} & {G_{T'}(f)}}\] commutes.
\end{prop}
We write $G(f)$ for $G_{T}(f)$, defined up to conjugation in $T$.
Abstractly, $G(f)$ is a quotient of $\pi_1(\M,*)$, a free group of
rank $|P|-1$.

\subsection{Recurrence}
Choose now for each $v\in T_1$ a path $\ell_v$ from $*$ to $v$ in
$\M$. Consider a loop $\gamma$ at $*$; then it induces a permutation
of $T_1$ given by $v\to v^\gamma$, and for each $v\in T_1$ its lift
$\gamma_v$ at $v$ yields a loop $\ell_v\gamma_v\ell_{v^\gamma}^{-1}$
at $*$, again written $\gamma_v$.
\begin{prop}
  $\gamma_v$ depends only on the class of $\gamma$ in $G(f)$.
\end{prop}

We therefore have a natural wreath product embedding
\[\phi:G(f)\to G(f)\wr\sym{T_1}\]
mapping $\gamma$ to $(v\to\gamma_v,v\to v^\gamma)$. Enumerating
$T_1=\{v_1,v_2,\dots,v_d\}$, we write
\[\phi(g)=\pair|g_{v_1},\dots,g_{v_d}|\pi_g,\]
where $\pi_g$ is the permutation of $T_1$ corresponding to $g$, in
disjoint cycle notation.

Since $G(f)$ is finitely generated (by at most $|P|-1$ elements), it
can be completely described by a permutation of $d$ points and a
sequence of $d$ elements of $G(f)$, written as products of generators.
We will call this description a \emph{wreath presentation}.

\subsection{Main example}
We study now the main example $f(z)=z^2-1$, and derive an explicit
action of $G=G(f)$ on a ``standard'' binary rooted tree $\{x,y\}^*$.

The postcritical set $P$ is $\{0,1,\infty\}$, so $\M$ is a
thrice-punctured sphere and $G$ is $2$-generated.

For convenience, pick as base point $*$ a point close to
$(1-\sqrt5)/2$; then $T_1=\{x,y\}$ with $x$ close to $*$ and $y$ close
to $-*$.

Consider the following representatives of $\pi_1(\M,*)$'s generators:
$a$ is a straight path approaching $-1$, turning a small loop in the
positive orientation around $-1$, and returning to $*$. similarly, $b$
is a straight path approaching $0$, turning around $0$ in the positive
orientation, and returning to $*$.

Let $\ell_x$ be a short arc from $*$ to $x$, and let $\ell_y$ be a
half-circle above the origin from $*$ to $y$.

Let us compute first $f^{-1}(a)$, i.e.\ the path traced by
$\pm\sqrt{z+1}$ as $z$ moves along $a$. Its lift at $x$ moves towards
$0$, passes below it, and continues towards $y$. Its lift at $y$ moves
towards $0$, passes above it, and continues towards $x$. We have
$a_x=b$ and $a_y=1$, so the wreath presentation of $a$ is
$\phi(a)=\pair|b,1|(x,y)$.

Consider next $f^{-1}(b)$. Its lift at $x$ moves towards $-1$, loops
around $-1$, and returns to $x$. Its lift at $y$ moves towards $1$,
loops and returns to $y$. We have $b_x=a$ and $b_y=1$, so the wreath
presentation of $b$ is $\phi(b)=\pair|a,1|$.

\section{Algebraic Properties}
We will exploit as much as possible the wreath presentation of
$G=G(z^2-1)$ to obtain algebraic information on $G$. Recall that $G$
is given with an action on the tree $T=\{x,y\}^*$. We write $xT,yT$
the subtrees below $x,y$ respectively.

Most of our proofs are close adaptations
of~\cite{brunner-s-v:nonsolvable}.

\begin{prop}\label{prop:wb}
  $G$ is weakly branch, i.e.\ in $G$'s action on $T$
  \begin{enumerate}
  \item $G$ acts transitively on each layer $T_n$;
  \item at every vertex $v$, there is a non-trivial $g\in G$ acting
    only on the subtree below $g$, and fixing its complement.
  \end{enumerate}
\end{prop}
\begin{proof}
  $G$ acts transitively on $T_1$, since $a$ permutes its two elements.
  Now $\stab(x)=\langle a^2,b,b^a\rangle$ maps onto $G$ by restriction
  of the action to $xT$; hence by induction acts transitively on each
  of its layers. Combining with the transitive action on $T_1$, we get
  a transitive action on all layers of $T$.
  
  Consider next $G'$. From $[a,b]=\pair|a^{-1}b,b^{-1}a|$ we get
  $G'\neq1$. From the relation $[b,a^2]=\pair|[a,b],1|$, we see that
  $G'$ contains $G'\times1$ as a subgroup acting only on $xT$. This
  can be iterated, so $G$ contains a copy of $G'$ acting only below
  the vertex $v$.
\end{proof}

Consider again the embedding $\phi:G\to G\wr\sym{x,y}$; define a norm
on the generators of $G$ by $\|a\|=\|a^{-1}\|=1$ and
$\|b\|=\|b^{-1}\|=\sqrt2$, and extend it by the triangle inequality:
\[\|g\|=\min\{\|s_1\|+\dots+\|s_n\|:\,g=s_1\dots
s_n,\,s_i\in\{a^{\pm1},b^{\pm1}\}\}.\]
\begin{lem}\label{lem:contr}
  If $\phi(g)=\pair|g_x,g_y|\pi$, then
  $\|g_i\|\le\frac{\|g\|+1}{\sqrt2}$ for all $i\in\{x,y\}$.
\end{lem}
\begin{proof}
  Write in a minimal form $g$ as a product of $N$ letters $a$ (with
  signs) and $M$ letters $b$ (with signs). Then $g_i$ contains at most
  $(N+1)/2$ letters $b$ and $M$ letters $a$.
\end{proof}

From now on write $c=[a,b]$.
\begin{prop}\label{prop:firstq}
  $G'=\langle c,c^a,c^{a^2}\rangle$; we have $G/G'=\Z^2$ generated by
  $a,b$ and $G'/(G'\times G')=\Z$ generated by $c$.
\end{prop}
\begin{proof}
  The first claim follows from computations: $c^b=c$, and
  $c^{a^3}=(c^{a^2})^{-1}cc^a$.
  
  Next, we have $G'=(G'\times G')\langle c\rangle$; indeed
  $c^a=\pair|[b,a],1|c^{-1}$, and
  $c^{a^2}=\pair|[b,a^{-1}],[b,a]^b|c$.
  
  Now assume for contradiction that $a^mb^n\in G'$ with $|m|+|n|$
  minimal. Then clearly $m$ is even, say $m=2p$. We have, for some
  $k\in\Z$ with $|k|\le|m|+|n|$,
  \[a^mb^n=\pair|b^pa^n,b^p|=\pair|g,h|c^k=\pair|ga^k,ha^{-bk}|\]
  and therefore $a^{n-k}b^p$ and $a^kb^p$ both belong to $G'$. This
  contradicts our assumption on minimality, so we have proven the
  second claim.

  Assume finally that $c^k\in G'\times G'$ with $|k|$ minimal. Then
  $a^k\in G'$ which contradicts the second claim.
\end{proof}

\begin{prop}\label{prop:tf}
  $G$ is a torsion-free group.
\end{prop}
\begin{proof}
  Since $G$ acts on the binary tree, it is residually a $2$-group, and
  its only torsion must be $2$-torsion. Assume for contraction that
  $G$ contains an element $g$ of order $2$, of minimal norm.
  
  By Proposition~\ref{prop:firstq} $a$ and $b$ are of infinite order.
  We may therefore assume $\|g\|\ge2$.  If $g$ fixes $x$, then its
  restrictions $g_x$ and $g_y$ are shorter by Lemma~\ref{lem:contr},
  and one of them has order $2$, contradicting $g$'s minimality.
  
  If $g$ does not fix $x$, then we may write $g=\pair|g_x,g_y|a$ for
  some $g_x,g_y\in G$. We then have $h=g^2=\pair|g_x bg_y,g_yg_xb|=1$
  and therefore $g_yg_xb=1$. Now for any element $h$ fixing $x$ we
  have $h_xh_y\in\langle a,b^2,G'\rangle$; this last subgroup does not
  contain $b$ by Proposition~\ref{prop:firstq}, so we have a contradiction.
\end{proof}

\begin{prop}\label{prop:fm}
  The subsemigroup $\{a,b\}^*$ of $G$ is free; therefore $G$ has
  exponential growth.
\end{prop}
\begin{proof}
  Consider two words $u,v$ in $\{a,b\}^*$ that are equal in $G$, and
  assume $|u|+|v|$ is minimal. We have $u_x=v_x$ and $u_y=v_y$ in $G$,
  which are shorter relations, so we may assume these words are equal
  by induction.
  
  Now if $u_x$ and $v_x$ start with the same letter, this implies that
  $u$ and $v$ also start with the same letter, and cancelling these
  letters would give a shorter pair of words $u,v$ equal in $G$.
\end{proof}

\begin{prop}\label{prop:lpres}
  $G$ has finite $L$-presentation
  \[G = \langle a,b|\,[[a^p,b^p],b^p],[[b^p,a^{2p}],a^{2p}]\text{
    for all $p$ a power of $2$}\rangle.\]

  The Schur multiplier $H_2(G,\Z)$ is free abelian of infinite rank.
\end{prop}
\begin{proof}
  This follows from the standard method in~\cite{bartholdi:lpres},
  since $G$ is contracting and $G'$ is finitely generated.
  
  We first obtain all relations $[[a^i,b],b]$ for odd $i$, and their
  iterates under the substitution $\sigma:a\mapsto b,b\mapsto a^2$
  induced by the inclusion of $G'\times1$ in $G'$.
  
  Then we note that $[[a^i,b],b]$ is a consequence of $[[a,b],b]$ and
  $[[b,a^2],a^2]$.
  
  The Schur multiplier computation follows from Hopf's formula
  $H_2(G,\Z)=(R\cap[F,F])/[R,F]$ for a group $G=F/R$ presented as
  quotient of a free group. The homomorphism $\sigma$ induces the
  standard shift on $\Z^\infty$.
\end{proof}

The method of the following proof is inspired by the PhD thesis of
Edmeia Da Silva~\cite{silva:phd}.
\begin{prop}\label{prop:nofg}
  $G$ does not contain any non-abelian free subgroup.
\end{prop}
\begin{proof}
  Take $g,h\in G$; we seek a relation $w(g,h)$ between them. Consider
  the standard metric $|g|$ on words, given by assigning length $1$ to
  each generator $a^{\pm1},b^{\pm1}$. Given any $g\in G$, we have
  $|g_x|\le|g|$ and $|g_y|\le|g|$; these inequalities are strict,
  unless $g\in B=\{a^{-i}b^na^j\text{ with }n\in\N,\,i,j\in\{0,1\}\}$.
  We also have $|g_xg_y|\le|g|$.
  
  We define the following ordering on pairs $(g,h)$ of group elements:
  Write $H=\stab_G(x)$ the stabilizer of $x$. The height of $(g,h)$ is
  $L(g,h)=|g|+|h|-\frac13\#\{g,h\}\cap H$; that is, we order first by
  total length, and then by number of elements in $\{g,h\}$ fixing
  $x$.
  
  The basis of the induction is that for $(g,h)$ of height at most $2$
  these elements must be generators, and we have the relation, say,
  $[[a,b],b]$ among them.
  
  If $g,h\in H$, they may be written as $g=\pair|g_x,g_y|$ and
  $h=\pair|h_x,h_y|$. Then either one of $g_x,h_x$ is not in $B$, so
  by induction there is a relation $u(g_x,h_x)$; or we consider the
  cases $g_x=b^i,h_x=b^j$ and $g_x=a^{-1}b^ia,h_x=a^{-1}b^ja$ which
  have a common power; and up to symmetry $g_x=a^{-1}b^ia,h_x=b^j$
  when we have the relation $[g_x,h_x]$.  By the same argument, there
  is a relation $v(g_y,h_y)$. Consider now relation $w=[u,v]$; we have
  $w(g,h)=[\pair|1,u(g_y,h_g)|,\pair|v(g_x,h_x),1|]=1$.
  
  If at least one, say $g$, fixes $x$, we may write
  $h=\pair|h_x,h_y|(x,y)$, and may proceed as above with $g$ and
  $h^2$, since $L(g_x,h_xh_y)<L(g,h)$ and $L(g_y,h_xh_y)<L(g,h)$.
  
  If none of $g,h$ fixes $x$, we write $g=\pair|g_x,g_y|(x,y)$ and
  $h=\pair|h_x,h_y|(x,y)$, and consider the elements in the following
  table:
  \[\begin{array}{cc|cc|cc|}
    &&\multicolumn{2}{c}{g_x\in H} &\multicolumn{2}{c}{g_x\not\in H}\\
    &&g_y\in H&g_y\not\in H&g_y\in H&g_y\not\in H\\ \hline
    &h_y\in H&g^2,h^2&g^2,h^2&g^2,h^2&g^2,h^2\\
    \raisebox{2ex}{$h_x\in H$\kern-1em}
    &h_y\not\in H&g^2,h^2&gh^{-1},hg^{-1}&gh,hg&g^2,h^2\\ \hline
    &h_y\in H&g^2,h^2&gh,hg&gh^{-1},hg^{-1}&g^2,h^2\\
    \raisebox{2ex}{$h_x\not\in H$\kern-1em}
    &h_y\not\in H&g^2,h^2&g^2,h^2&g^2,h^2&g^2,h^2\\ \hline
  \end{array}\]
  By inspection, in the pair of elements $e,f$ corresponding to $g,h$
  at both fix $x$, and their respective projections satisfy
  $L(e_x,f_x)<L(g,h)$ and $L(e_y,f_y)<L(g,h)$; therefore they satisfy
  relations $u,v$, and $e,f$ satisfy the relation $[u,v]$.
\end{proof}

\begin{prop}\label{prop:wcp}
  The word problem is solvable in $G$.
\end{prop}
\begin{proof}
  The solution to the word problem follows from Lemma~\ref{lem:contr};
  indeed take a word $w$ of length $n$. If it does not fix $x$, then
  it is not trivial. Otherwise its projections $w_x,w_y$ are shorter
  words, so can inductively be checked for triviality in $G$. Then $w$
  is trivial in $G$ if and only if $w_x$ and $w_y$ are both trivial in
  $G$.
\end{proof}

Write now $d=[c,a]$ and $e=[d,a]$.
\begin{prop}\label{prop:lcs}
  We have $G''=(\gamma_3\times\gamma_3)$, and
  $\gamma_3=(\gamma_3\times\gamma_3)\langle d,e rangle$.
  
  In the lower central series of $G$ we have $\gamma_1/\gamma_2=\Z^2$,
  $\gamma_2/\gamma_3=\Z$, and $\gamma_3/\gamma_4=\Z/4$. Therefore all
  quotients except the first two in the lower central series are
  finite.
  
  In the lower $2$-central series defined by $G_1=G$ and
  $G_{n+1}=[G,G_n]G_{n/2}^2$, we have
  \[\dim_{\F_2}\Gamma_n/\Gamma_{n+1}=\begin{cases}i+2& \text{ if
      }n=2^i\text{ for some }i;\\ 
    \max\setsuch{i+1}{2^i\text{ divides }n}& \text{ otherwise}.
  \end{cases}\]
\end{prop}
\begin{proof}
  We moreover claim we have the following generating sets:
  \begin{align*}
    \gamma_1 &= \langle a,\,b\rangle;\\
    \gamma_2 &= \langle c=[a,b]=(a,a^{-b}),\,
    c^{-1-a}=(c,1),\,c^{-a^{-1}-1}=(1,c)\rangle;\\
    \gamma_3 &= \langle d=[c,a],\,e=[d,a],\,[e^{-1},b]=(d,1),\,(e,1),\,(1,d),\,(1,e)\rangle;\\
    \gamma_4 &= \langle d^4,\,e,(d,1),\,(e,1),\,(1,d),\,(1,e)\rangle.
  \end{align*}
  Most follow from simple computations; let us justify
  $d^4\in\gamma_4$. Since it is clear that $\gamma_4$ contains $e$, we
  denote by $\equiv$ congruence modulo $e$, which is weaker than
  congruence modulo $\gamma_4$. We have
  \begin{align*}
    d^2 &\equiv d^2e=b^{-1}a^{-1}ba^{-2}b^{-1}aba^2\\
    &=b^{-1}aba^{-2}b^{-1}a^{-1}ba^2\text{ modulo }[a^{2b},a^2]=1\\
    &=(d^2e)^{-a^{-b}}\equiv d^{-2}.
  \end{align*}
  
  We have already computed $G/G'=\Z^2$ in
  Proposition~\ref{prop:firstq}, and $G'/\gamma_3=\Z$ is generated by
  $c$. This last computation gives $\gamma_3/\gamma_4=\Z/4$.

  Next, $G''$ is generated as a normal subgroup by $[c,(c,1)]=(d,1)$,
  so $G''=\gamma_3\times\gamma_3$.

  For the $2$-central series see~\cite{bartholdi:bsvlcs}, where the
  argument is developed for the BSV group.
\end{proof}

\begin{prop}\label{prop:jns}
  $G$ is not solvable, but all its quotients are (free
  abelian)-by-(finite $2$-group).
\end{prop}
\begin{proof}
  By Proposition~\ref{prop:lcs}, we have $G''>\gamma_3\times\gamma_3$,
  so $G'''>G''\times G''$, and hence $G^{(n)}>G^{(n-1)}\times
  G^{(n-1)}$ for all $n$. Assume for contradiction that $G$ is
  solvable; this means $G^{(n)}=1$ for some minimal $n$, a
  contradiction with the above statement.
  
  Now consider a non-trivial normal subgroup $N$ of $G$.
  By~\cite[Theorem~4]{grigorchuk:jibg}, we have
  $(\gamma_3)^{\times2^n}<N$ for some $n$. On the other hand, $G$ is
  an abelian-by-(finite $2$) extension of $(\gamma_2)^{\times2^n}$.
  The result follows.
\end{proof}

\section{Geometric and Analytic Properties}
In this part we isolate results on the geometry of the action of $G$
on the boundary of the tree. The inspiration for this part is the
paper~\cite{bartholdi-g:spectrum} where computations of spectra were
performed, and~\cite{bartholdi-g-n:fractal} where the relation
between the Schreier graphs of an i.m.g.\ group and the Julia set of
its polynomial were explicited.

Consider first the spectrum computation. Given a group $G$ with
generating set $S$ and a unitary representation $\pi$ on a Hilbert
space $\hilb$, its \emph{spectrum} is
\[\spec(\pi)=\frac1{|S|}\spec_{B(\hilb)}\left(\sum_{s\in
    S}\pi(s)\right).
\]
It is a closed subset of $\C$; if $S$ is symmetric, then the spectrum
is always contained in $\R$.

The main representation of interest is the left-regular representation
$\rho$ on $\ell^2(G)$. A result of Kesten~\cite{kesten:amen} shows
that $1$ is in $\spec(\rho)$ if and only if $G$ is amenable.

Here we concentrate on the ``natural'' representation $\pi$ of $G$ on
$L^2(\partial T$, where $\partial T$ is the boundary of the tree $T$.
We denote by $\pi_n$ the representation of $G$ on the $n$th layer
$T_n$.

Introduce for $n\ge0$ the following homogeneous polynomials:
\[Q_n(\lambda,\mu,\nu)=\det\left(\lambda+\mu(\pi_na+\pi_na^{-1})+\nu(\pi_nb+\pi_nb^{-1})\right).\]
Then the spectrum of $\pi_n$ is obtained by solving
$Q_n(\lambda,-\frac14,-\frac14)=0$.

\begin{figure}
  \label{fig:spectrum}
  \psfrag{\2611}{$-1$}
  \psfrag{\2610.8}{$0.8$}
  \psfrag{\2610.6}{$0.6$}
  \psfrag{\2610.4}{$0.4$}
  \psfrag{\2610.2}{$0.2$}
  \psfrag{0}{$0$}
  \psfrag{0.2}{$0.2$}
  \psfrag{0.4}{$0.4$}
  \psfrag{0.6}{$0.6$}
  \psfrag{0.8}{$0.8$}
  \psfrag{1}{$1$}
  \begin{center}
    \epsfig{file=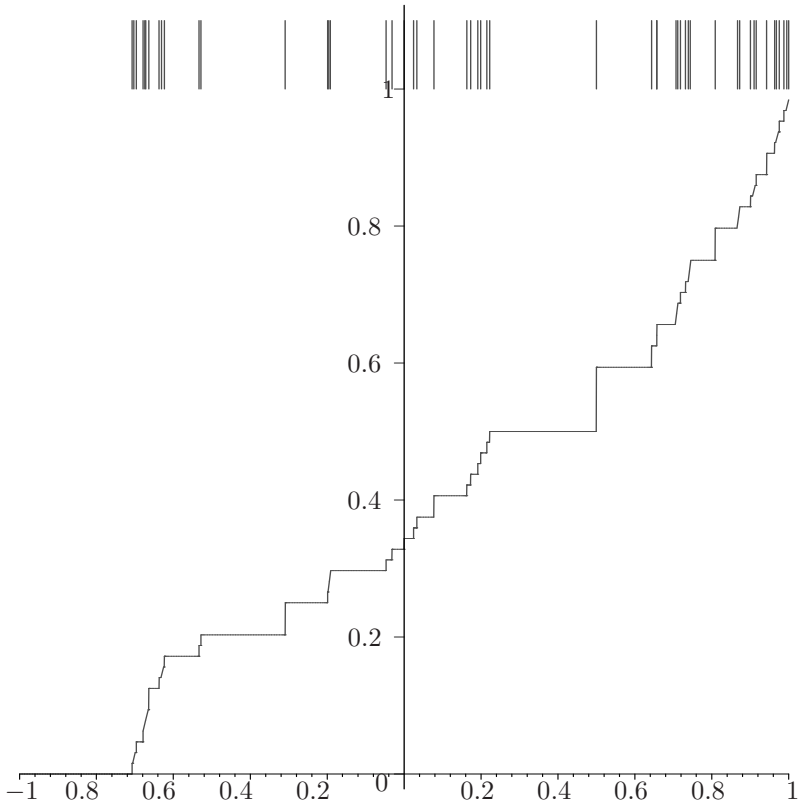,width=300pt,height=300pt}
  \end{center}
  \caption{The spectrum and spectral measure of $\pi$, in its
    level-$6$ approximation.}
\end{figure}

We are not able to derive a closed form for $Q_n$; however we have the
\begin{prop}\label{prop:spectrum}
  Define the polynomial mapping $F:\R^3\to\R^3$ by
  \[(\lambda,\mu,\nu)\mapsto(\lambda^2+2\lambda\nu-2\mu^2,\lambda\nu+2\nu^2),-\mu^2).\]
  Then $Q_n$ is given by
  \begin{gather*}
    Q_0(\lambda,\mu,\nu)=\lambda+2\mu+2\nu;\\
    Q_1(\lambda,\mu,\nu)=Q_0(\lambda,\mu,\nu)\cdot(\lambda-2\mu+2\nu);\\
    \intertext{and for $n\ge1$,}
    Q_{n+1}(\lambda,\mu,\nu)=Q_n(F(\lambda,\mu,\nu)).
  \end{gather*}
  
  Define $K$ as the closure of the set of all backwards $F$-iterates
  of $\{Q_1=0\}$. Then the spectrum of $\pi$ is the intersection of
  $\{\mu=\nu=-\frac14\}$ with $K$.
\end{prop}
\begin{proof}
  The first two follow from direct computation, using
  $\pi_0(a)=\pi_0(b)=(1)$, and
  \[\pi_{n+1}(a)=\begin{pmatrix}0&\pi_n(b)\\1&0\end{pmatrix},\quad
  \pi_{n+1}(b)=\begin{pmatrix}\pi_n(a)&0\\0&1\end{pmatrix}.\]
  Next, we compute, writing $a,b$ for $\pi_n(a),\pi_n(b)$ respectively,
  \begin{align*}
    Q_{n+1}(\lambda,\mu,\nu)&=\det\begin{pmatrix}\lambda+\nu(a+a^{-1})
      & \mu(1+b)\\ \mu(b^{-1}+1) & \lambda+2\nu\end{pmatrix}\\
    &=
    \det\left((\lambda+\nu(a+a^{-1}))(\lambda+2\nu)-\mu^2(1+b)(b^{-1}+1)\right)\\
    &= Q_n(\lambda^2+2\lambda\nu-2\mu^2,\lambda\nu+2\nu^2),-\mu^2).
  \end{align*}
  The claim follows.
\end{proof}
It seems difficult to obtain a more explicit description of
$\spec(\pi)$, since the polynomials $Q_n$ do not factor well. For
instance, $Q_5$ has an irreducible factor of degree $5$, and
$Q_5(\lambda,-\frac14,-\frac14)$ has Galois group $S_5$.
Figure~\ref{fig:spectrum} displays the spectrum of $\pi_6$ and its
associated spectral measure; from there we see that the spectrum of
$\pi$ is a Cantor set.

Now we turn our attention to the Schreier graphs of $G$. These are
graphs $\gf_n$ with $2^n$ points corresponding to the vertices in
$T_n$, with edges labelled $a$ and $b$ connecting $v$ to $v^a,v^b$
respectively for all $v\in T_n$.

These graphs are therefore $4$-regular graphs; the Green function of
the random walk on $\gf_n$ is given by the integral of $1/(1-t)$ with
respect to the spectral measure given above. The striking fact is the
following:
\begin{prop}\label{prop:limit}
  The graphs $\gf_n$ are planar, and can be drawn in $\C$ in such a
  way that they converge in the Hausdorff metric to the Julia set $\J$
  of $z^2-1$.
\end{prop}
\begin{proof}
  This follows from~\cite{bartholdi-g-n:fractal}, since $z^2-1$ is
  hyperbolic.
\end{proof}

$\gf_n$ is constructed as follows: it is built of two parts $A_n,B_n$
connected at a distinguished vertex. Each of these is $4$-regular,
except at the connection vertex where each is $2$-regular, and $A_n$
contains only the $a^{\pm1}$-edges while $B_n$ contains only the
$b^{\pm1}$-edges.

$A_0$ and $B_0$ are the graphs on $1$ vertex with a single loop of the
appropriate label.
  
If $n=2k$ is even, then $B_{2k+1}=B_{2k}$, and $A_{2k+1}$ is obtained
by taking an $a$-labelled $2^{k+1}$-gon $v_0,\dots,v_{2^{k+1}-1}$, and
attaching to each $v_i$ with $i\neq0$ a copy of $B_{2j}$ where $j$ is
the largest power of $2$ dividing $i$. Its distinguished vertex is
$v_0$.

If $n=2k-1$ is odd, then $A_{2k}=A_{2k-1}$, and $B_{2k}$ is obtained
by taking an $b$-labelled $2^k$-gon $v_0,\dots,v_{2^k-1}$, and
attaching to each $v_i$ with $i\neq0$ a copy of $A_{2j+1}$ where $j$
is the largest power of $2$ dividing $i$. Its distinguished vertex is
$v_0$.

\begin{figure}
  \label{fig:schreier}
  \begin{center}
    \unitlength=1pt\begin{picture}(0,0)
      \put(15,131){{\small $g_1$}}
      \put(35,131){{\small $g_2$}}
      \put(57,110){{\small $g_1$}}
      \put(57,142){{\small $g_1$}}
      \put(94,110){{\small $g_1$}}
      \put(94,142){{\small $g_1$}}
      \put(121,131){{\small $g_2$}}
      \put(136,131){{\small $g_1$}}
      \put(36,103){$\Gamma_3$}
      \put(223,103){$\Gamma_4$}
      \put(45,16){$\Gamma_5$}
      \put(220,16){$\Gamma_6$}
    \end{picture}
    \epsfig{file=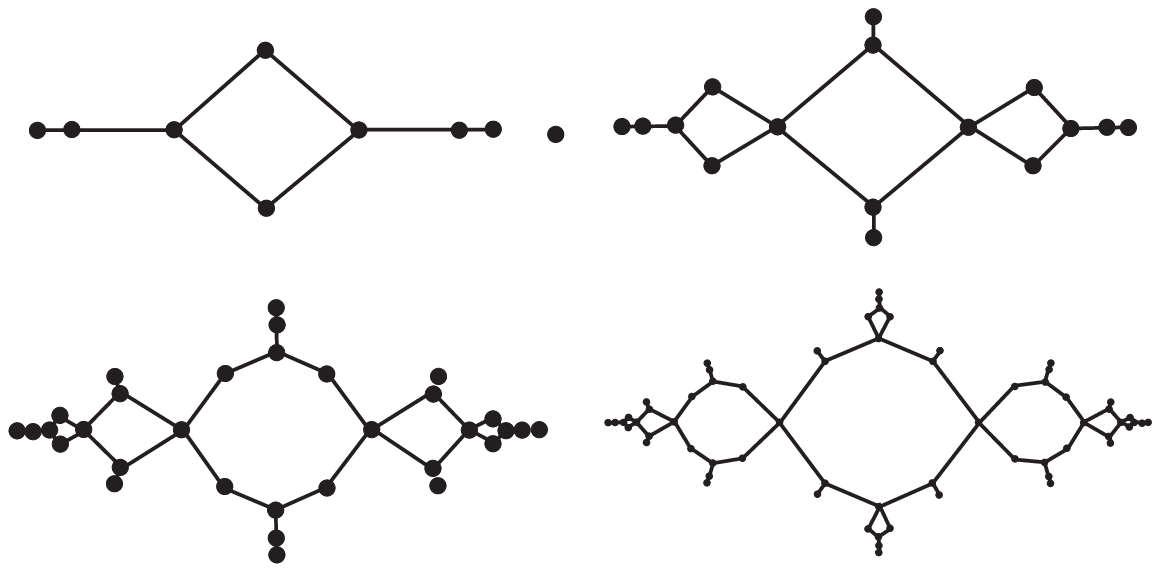}
  \end{center}
  \caption{The Schreier graphs $\gf_n$ for $3\le n\le 6$.}
\end{figure}

\begin{figure}
  \label{fig:julia}
  \begin{center}
    \epsfig{file=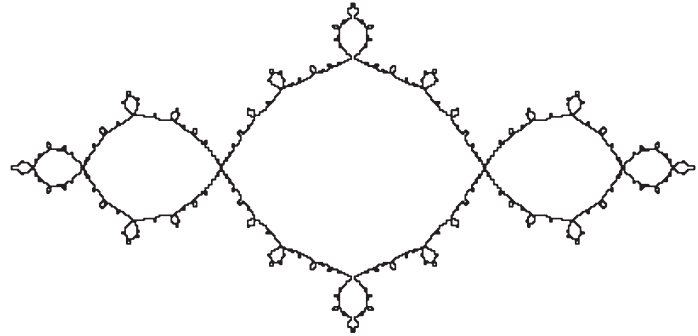,width=100mm}
  \end{center}
  \caption{The Julia set $\J$ of the polynomial $z^2-1$.}
\end{figure}

The first Schreier graphs $\gf_n$ of $G$ are drawn in
Figure~\ref{fig:schreier}. Compare with the Julia set in
Figure~\ref{fig:julia}.

\bibliography{mrabbrev,people,math,grigorchuk,bartholdi}
\end{document}